\documentclass[12pt,reqno,twoside]{article}
\usepackage{amsfonts,amssymb,latexsym,xspace,epsfig,graphics,color}
\usepackage{amsmath,amsthm,enumerate,stmaryrd,xy}
\usepackage{bbm}
\usepackage{natbib}
\usepackage{subfigure}
\usepackage{times}

\def\titlerunning#1{\gdef\titrun{#1}}
\makeatletter
\def\author#1{\gdef\autrun{\def\and{\unskip, }#1}\gdef\@author{#1}}
\def\address#1{{\def\and{\\\hspace*{18pt}}\renewcommand{\thefootnote}{}%
\footnote {#1}}%
\markboth{\autrun}{\titrun}}
\makeatother
\def\email#1{e-mail: #1}
\def\subjclass#1{{\renewcommand{\thefootnote}{}%
\footnote{\emph{Mathematics Subject Classification (2010):} #1}}}
\def\keywords#1{\par\medskip
\noindent\textbf{Keywords.} #1}

\newcommand{\x}{\ensuremath{\underline{x}}}
\newcommand{\y}{\ensuremath{\underline{y}}}

\newcommand{\F}{\mathcal{F}}

\newcommand{\sig}{\sigma}
\newcommand{\N}{\ensuremath{\mathbb{N}}}

\newcommand{\R}{\ensuremath{\mathbb{R}}}     
\newcommand{\Z}{\ensuremath{\mathbb{Z}}}

\def\eqd{\,{\buildrel D \over =}\,} 

\newcommand{\X}{\ensuremath{\mathcal{X}}}

\def\vareps{\varepsilon}

\usepackage{algorithmic}
\usepackage{algorithm}

\newtheorem{lemma}{Lemma}

\newtheorem*{defi}{Definition}
\newtheorem{theo}{Theorem}
\newtheorem*{theo*}{Theorem}

\newtheorem{coro}{Corollary}

\begin{document}

\baselineskip=17pt

\titlerunning{Characterization of the stability of chains associated with $g$-measures}

\title{Characterization of the stability of chains associated with $g$-measures}

\author{Christophe Gallesco, Sandro Gallo $\&$ Daniel Y. Takahashi}

\date{}
\maketitle

\address{Departmento de Estat\'istica, Instituto de Matem\'atica, Estat\'istica e Ci\^encia de Computa\c{c}\~ao, Universidade de Campinas, Brasil; \email{gallesco@ime.unicamp.br}; \and Departamento de Estat\'istica, Universidade Federal de S\~ao Carlos, Brasil; \email{sandrodobrasil@gmail.com} \and Neuroscience Institute and Psychology Department, Princeton University, USA; \email{takahashiyd@gmail.com}}



\begin{abstract}

In this paper we introduce a notion of asymptotic stability of a probability kernel, which we call dynamic uniqueness.
We say that a kernel exhibits dynamic uniqueness if all the stochastic chains starting from a fixed past coincide on the future tail $\sigma$-algebra. We prove that the dynamic uniqueness is generally stronger than the usual notion of uniqueness for $g$-measures. Our main result shows that dynamic uniqueness is equivalent to the weak-$\ell^2$ summability condition on the kernel. This generalizes and strengthens the Johansson-\"Oberg $\ell^2$ criterion for uniqueness of $g$-measures. Finally, among other things, we prove that the weak-$\ell^2$ criterion implies $\beta$-mixing of the unique $g$-measure compatible with a regular kernel improving several results in the literature.  


\end{abstract}


\section{Introduction}
 
Let $S$ be a finite set, $\mathcal{X} = S^{\Z}$ and  $\mathcal{X}^- = S^{\Z_-^*}$, where $\Z_-^* = \{-1, -2, \ldots\}$. We denote by $x_i$ the $i$-th coordinate of $x \in \mathcal{X}$ and for $-\infty \leq j \leq i < \infty$ we write $x^{i}_{j}:=(x_{i},\ldots, x_{j})$ and $x^{i}_{-\infty} :=(x_{i},x_{i-1},\dots)$. We also use the shorthand notation $\x:=x_{-\infty}^{-1}$. For $x, y \in \mathcal{X}$ and $i, j,k$ finite, $j\leq i$, a \emph{concatenation} $y^{i}_{j}x^k_{-\infty}$ is a new sequence $z\in S^{\{ i, i-1, \ldots\}}$ with $z^{i}_{j} = y^{i}_{j}$ and $z^{j-1}_{-\infty} = x^k_{-\infty}$. 
Note that we are using the convention that the past (smaller indices) of $x \in \mathcal{X}$ is represented on the right hand side.

For each $\Lambda \subset \Z$, consider the set $\mathcal{X}_\Lambda = S^{\Lambda}$, together with the canonical projection $\pi_\Lambda:  \mathcal{X}\rightarrow \mathcal{X}_\Lambda$ defined by $\pi_\Lambda(x)_k = x_k$ for all $k \in \Lambda$. Then, define
$\mathcal{C}(\Lambda) = \{ \pi^{-1}(B): B \subset\mathcal{X}_\Lambda\}$, called the cylinders with base $\Lambda$. For $\Gamma \subset \Z$, we consider the algebra of cylinders with base in $\Gamma$ defined by $\mathcal{C}_\Gamma = \bigcup\{\mathcal{C}(\Lambda): \Lambda \subset \Gamma, \Lambda \;\text{finite}\}$ and the $\sigma$-algebra generated by the algebra of cylinders with base in $\Gamma$, $\F_\Gamma = \sigma(\mathcal{C}_\Gamma)$. We use the shorthand notation $\F = \F_\Z$,  $\F^- = \F_{\Z_-^*}$, $\F^+ = \F_{\Z_+}$, and $\F_{I} = \F_{I \cap \Z}$, for any $I \subset \R$. Let $\mathcal{T}^+ = \bigcap_{n \geq 1} \F_{[n, \infty)}$ be the future tail $\sigma$-algebra.

A \emph{probability kernel}, or simply a \emph{kernel}  $g$ on the alphabet $S$ is a measurable function  
\begin{equation*}
\begin{array}{cccc}
g:&\mathcal{X}^-&\rightarrow& [0,1]\\
\end{array}
\end{equation*}
such that
\[
\sum_{a\in S}g(a\x)=1\,\,,\,\,\,\,\,\,\forall \x \in \mathcal{X}^-.
\] 
		
We say that a stationary stochastic process $(X_{j})_{j\in\Z}$ with values in $S$ defined on a probability space $(\mathcal{X}, \mathcal{F},P)$  is \emph{compatible} with a kernel $g$ if the latter is a regular version of the conditional probabilities of the former, \textit{i.e.},
\begin{equation*}\label{compa}
P[X_{0}=a\mid X_{-\infty}^{-1}=\x]=g(a\x)
\end{equation*}
for every $a\in S$ and $P$-a.e.\ $x$ in $\mathcal{X}$. The measure $P$ is called $g$-measure. 
A kernel $g$ is  \emph{strongly non-null} if 
\begin{equation}
\label{strongnonnull}
\inf_{a \in S,\x\in \mathcal{X}^-}g(a\x) > 0.
\end{equation}
The \emph{variation rate} (or \emph{continuity rate}) \emph{of order} $k$ of a kernel $g$ is given by
\begin{equation} 
\label{eq:Vardef}
\textrm{var}_{k}(g):=\sup_{a\in S}\sup_{\omega_{-k}^{-1}\in S^{k}}\sup_{\x,\y\in \mathcal{X}^-}|g(a\omega_{-k}^{-1}\x)-g(a\omega_{-k}^{-1}\y)|. 
\end{equation}
We say that $g$ is continuous if $\lim_{k\rightarrow \infty}\textrm{var}_{k}(g)=0$. A kernel $g$ is in $\ell^q$ if 
$$\sum_{k=1}^\infty \textrm{var}_{k}(g)^q < \infty.$$ 

If the kernel is continuous,  a compactness argument shows that we can take the weak limits of the processes starting at fixed pasts (see below) to obtain compatible stationary chains (see for example \cite{keane/1972}). If $g$ is strongly non-null and continuous, we say that $g$ is a \emph{regular kernel}. 
When there is more than one stationary process compatible with $g$, we say that there is an \emph{equilibrium phase transition}, otherwise we say that the process is \emph{unique} (or that there is an \emph{equilibrium uniqueness}). This is the usual notion of uniqueness/phase transition of $g$-measures in the literature. 
		
Let $\mathcal{M}$ be the set of probability measures on $(\mathcal{X}^-, \mathcal{F}^-)$. Given a probability kernel $g$, we call $g$-chain with initial distribution $\lambda \in \mathcal{M}$, the canonical process $(X_{j})_{j\in\Z}$ on the probability space $(\mathcal{X}, \mathcal{F},P^\lambda)$  defined in the following way. For any event $B \in \mathcal{F}^-$, $P^\lambda[B] = \lambda[B]$ and,  for $j \geq 0$,
\begin{equation*}\label{compa}
P^{\lambda}[X_{j}=a\mid X_{-\infty}^{j-1}=y^{j-1}_{-\infty}]=g(ay^{j-1}_{-\infty}).
\end{equation*}
By Ionescu-Tulcea extension theorem, $P^\lambda$ is uniquely defined. The set of measures that can be written as $P^\lambda$ for some $\lambda \in \mathcal{M}$ is denoted $\mathcal{M}_0$. When  $\lambda[X^{-1}_{-\infty} = \x] = 1$ for some $\x \in \mathcal{X}^{-1}_{-\infty}$, we use the notation $P^{\x}$ for the associated probability measure. By standard arguments in measure theory we have that $P^\lambda[B] = \int_{\X}P^{\x}[B]\lambda[d\x]$ for all $B \in \F$. 

The purpose of this article is to study a notion of asymptotic stability for $g$-chains, which we call \emph{dynamic uniqueness}. The basic idea is the following. Given a kernel $g$ and the corresponding $g$-chains with distinct fixed pasts, we want to study their measures on the future tail $\sigma$-algebra $\mathcal{T}^+$. The events on this tail $\sigma$-algebra (tail events) can have the  physical interpretation of asymptotic events and are good candidates to detect phase transitions. We say that there is dynamic uniqueness if the measures on the tail $\sigma$-algebra agree for all pasts. The kernel is therefore asymptotically stable because, on the future tail $\sigma$-algebra, the specified chains do not depend on the initial condition. Conversely, we say that a kernel exhibits a dynamic phase transition if there exist two different pasts for which the corresponding measures on the tail events disagree.

\begin{defi}\label{def:dynamic}
We say that a kernel $g$ exhibits \emph{dynamic uniqueness} if for all events  $A \in \mathcal{T}^+$ and any pair $\x,\y \in \X^{-}$, we have $P^{\x}[A] = P^{\y}[A]$, otherwise we say that there is \emph{dynamic phase transition}. 
\end{defi}

We first prove several equivalent criteria for dynamic uniqueness. The main result is the equivalence between dynamic uniqueness and the weak-$\ell^2$ criterion. We note that Johansson and \"Oberg $\ell^2$ criterion \citep{johansson/oberg/2003} implies the weak-$\ell^2$ criterion. We also show that dynamic uniqueness is equivalent to convergence in total variation distance. This is akin to the convergence of Markov chains to the invariant distribution in total variation distance.  \cite{johansson/oberg/pollicott/2012} showed that the $\ell^2$ criterion implies the weak convergence of the $g$-chains to the unique equilibrium measure. Using our result, we can strengthen the weak convergence to the convergence in total variation distance.

Furthermore, using our criteria for dynamic uniqueness, we demonstrate that for the  Bramson-Kalikow-Friedli (BKF) \citep{bramson/kalikow/1993, friedli/2010} and  Hulse  \citep{hulse/2006} models, there exists dynamic uniqueness if and only if the rate of variation is in $\ell^{2}$. Because there are BKF and Hulse models that exhibit equilibrium uniqueness but do not have rate of variation in $\ell^{2}$, we conclude that equilibrium and dynamic uniqueness/phase transition can differ in general. 
Finally, we show that dynamic uniqueness implies that the $g$-measure has a strong form of mixing called weak Bernoullicity (which is equivalent to absolute regularity and also $\beta$-mixing). For the proof of this result, we exhibit several characterizations of $g$-measures that are weak Bernoulli. In particular, we show that a $g$-measure is weak Bernoulli if and only if the weak-$\ell^2$ criterion holds almost surely. This strengthens a result by \citet{walters/2005} that shows that a $g$-measure is weakly Bernoulli if the kernel is in the Bowen class ($Bow(\X^-, T)$).  Our result also strengthens a result by \cite{johansson/oberg/pollicott/2012} that shows that $\ell^2$ criterion implies that the unique $g$-measure is very weak Bernoulli. To conclude, we give an example of a kernel exhibiting dynamic phase transition but the corresponding $g$-measure is weak Bernoulli. This shows that, in general, dynamic uniqueness is a stronger condition than weak Bernoullicity.

\section{Results}
 In what follows, we collect some definitions that we will need to state the results.
A probability measure $\mu$ is \emph{trivial} on a $\sigma$-algebra $\mathcal{H}$, if for all $A \in \mathcal{H}$ we have $\mu(A) = 1$ or $0$. Let $\mathcal{H}'$ be a sub-$\sigma$-algebra of $\mathcal{H}$. We indicate by $\mu|_{\mathcal{H}'}$ the restriction of $\mu$ to $\mathcal{H}'$. Let $\tilde{\mu}$ be a coupling between $\mu$ and $\mu'$ and $(X_j, X'_j)_{j \in \Z}$ be the associated process, where $(X_j)_{j \in \Z} \eqd \mu$ and $(X'_j)_{j \in \Z} \eqd \mu'$ (the symbol ``$\eqd$" stands for equality in law). The \emph{coupling time} $\Theta$ is defined as $\Theta = \inf \{t\geq 0: X_j = X'_j\;\; \text{for all}\;\; j\geq t\}$. We say that a measure $\mu$ is \emph{mixing} if, for all $B \in \F$,
\begin{equation} \label{eq:mixing}
 \lim_{n \rightarrow \infty}\sup_{A\in \F_{[n, \infty)}}|\mu[A \cap B] - \mu[A]\mu[B] | = 0.
\end{equation}

\begin{defi}
We say that a regular kernel $g$ with a $g$-measure $P$  satisfies 
\begin{itemize}
\item the \emph{$\ell^2$ criterion} if the variation rate of $g$ is in $\ell^2$ or, equivalently,\begin{equation*}
\sum_{n=0}^\infty \sum_{a \in S}\sup_{\omega^{n}_{0} \in S^n}\sup_{\x, \y \in \X^-}\left(g(a\omega^{n}_0\x)-g(a\omega^{n}_0\y)\right)^2 < \infty;
\end{equation*}

\item the \emph{weak-$\ell^2$ criterion} if, for every pair $\x,\y \in \X^-$, 
\begin{equation*}
\sum_{n=0}^\infty \sum_{a \in S}\left(g(a\omega^{n}_0\x)-g(a\omega^{n}_0\y)\right)^2 < \infty, \;\textrm{for}\;P^{\x}\textrm{-a.e.}\,\, \omega;
\end{equation*}
\item the $P$-weak-$\ell^2$ criterion if the above inequality holds for $P|_{\F^-} \otimes P|_{\F^-}$-a.e.~$(\x,\y)$. 
\end{itemize}

\end{defi}

The $\ell^2$-criterion is well-known in the literature to be the tightest criterion for equilibrium uniqueness \citep{johansson/oberg/2003, johansson/oberg/pollicott/2012}. One of the main objectives of the present paper is to explore the two other criteria and their implications on the statistical properties of the associated measures. In case the reader wants to have a bird's-eye view of the main results, we summarized them in the end of this section. 

Our first result is the following theorem that characterizes dynamic uniqueness.
\begin{theo} \label{theo:dynamicunique}
Let $g$ be a regular kernel. The following statements are equivalent:
\begin{enumerate}[(i)]
 \item $g$ exhibits dynamic uniqueness. \label{1}
 
 \item For all $\x,\y \in \mathcal{X}^-$, \label{2}
$
 \lim_{n \rightarrow \infty} \sup_{B \in \F_{[n, \infty)}} |P^{\x}[B] - P^{\y}[B]| = 0.
$

\item For all $\x,\y \in \mathcal{X}^-$ there exist a coupling $\tilde{P}$ between $P^{\x}$ and $P^{\y}$ and a coupling time $\Theta$ such that $\tilde{P}[\Theta < \infty] = 1$. \label{3}

 \item For all $\lambda,\lambda' \in \mathcal{M}$, $P^{\lambda}|_{\mathcal{T}^+} = P^{\lambda'}|_{\mathcal{T}^+} $. \label{4}

\item For all $\lambda \in \mathcal{M}$, $P^\lambda$ is trivial on $\mathcal{T}^+$.   \label{5}

\item For all $\lambda \in \mathcal{M}$, $P^\lambda$ is mixing. \label{6}

\item For all $\x,\y \in \mathcal{X}^-$, $P^{\x}|_{\F^+} \ll P^{\y}|_{\F^+}$. \label{7}

\item  $g$ satisfies the weak-$\ell^2$ criterion. \label{8}

\end{enumerate}
\end{theo}

Equivalences $(\ref{1})$,$(\ref{2})$, and $(\ref{3})$ in Theorem \ref{theo:dynamicunique} follow from the characterization of exact maximal coupling given by \cite{thorisson/2000}.  We remind the reader that given two probability measures $\mu$ and $\nu$ on $(\X, \F)$, the total variation distance between $\mu$ and $\nu$ is given by 
$$d_{TV}(\mu, \nu) := 2\sup_{A \subset \F} |\mu(A) - \nu(A)|.$$  
Thus $(\ref{2})$ is equivalent to
$$\lim_{n \rightarrow \infty}d_{TV}(P^{\x}|_{\F_{[n, \infty)}}, P^{\y}|_{\F_{[n, \infty)}}) = 0$$ 
and demonstrates the convergence in total variation.  Equivalence $(\ref{4})$ shows that the probabilities of the tail events of the non-stationary and stationary $g$-chains coincide when there is dynamic uniqueness. Equivalences $(\ref{5})$ and $(\ref{6})$ show that dynamic phase transition is observed when there exists an initial condition for which the  correlations do not decay to zero, formalizing the idea that a phase transition happens when there are long range correlations. 

The reader familiar with the general theory of Markov chains will notice that equivalences $(i)$ to $(vi)$ were inspired by the characterization of total variation convergence for Markov chains. The main difference for the proof of these equivalences for $g$-chains is that we do not have the Markov property that is crucial to demonstrate the existence of a Markov chain with trivial future tail $\sigma$-algebra. We bypass the lack of Markovianess by using the Rohlin and Sinai Theorem \citep{rohlin/sinai/1961} to show that the past and future tail $\sigma$-algebras coincide almost surely for the $g$-measures. 
Theorem \ref{theo:dynamicunique} $(\ref{7})$ and $(\ref{8})$  are obtained using the theory of absolute continuity and singularity developed by Shiryaev and co-authors \citep{jacod/shiryaev/2002, engelbert/shiryaev/1980}. Condition $(\ref{8})$ is particularly useful as can be used to generate criteria for dynamic uniqueness as it is shown in Corollary \ref{theo:dynamicuniquel2} below. 

A corollary of the equivalences $(\ref{1})$ and $(\ref{2})$ in Theorem \ref{theo:dynamicunique}  is the following relation between dynamic and equilibrium phase transitions. 

\begin{coro} \label{coro:dynamicphase}
 If a kernel $g$ is regular and exhibits equilibrium phase transition, then it exhibits dynamic phase transition.
\end{coro}

\citet{johansson/oberg/2003} proved that if a regular kernel $g$ is in $\ell^2$, then it has a unique compatible $g$-measure. A consequence of Theorem \ref{theo:dynamicunique} $(\ref{8})$ is that the same criterion guarantees dynamic uniqueness.

\begin{coro} \label{theo:dynamicuniquel2}
If $g$ is a regular kernel that satisfies the $\ell^{2}$ criterion, then there is dynamic uniqueness.
\end{coro}

 \cite{johansson/oberg/pollicott/2012} proved that if $g$ is a kernel in $\ell^{2}$, then the $g$-chain starting with any initial measure converges in $\bar{d}$ distance to the unique compatible stationary measure. The $\bar{d}$ convergence is stronger than the weak convergence but strictly weaker than the total variation convergence. Therefore,  Corollary \ref{theo:dynamicuniquel2} together with Theorem \ref{theo:dynamicunique} $(ii)$ imply that we can strengthen the result in \citet{johansson/oberg/pollicott/2012} and obtain a convergence in total variation distance. Also, Corollary \ref{theo:dynamicuniquel2} implies that if a kernel is in $\ell^2$ then there is a coupling satisfying $(\ref{3})$  in Theorem \ref{theo:dynamicunique}. This is rather surprising given that couplings achieving $(\ref{3})$ were previously obtained only under summable continuity rate \citep{gallo/lerasle/takahashi/2013, comets/fernandez/ferrari/2002}. 

Corollaries \ref{coro:dynamicphase} and \ref{theo:dynamicuniquel2} together imply that models exhibiting equilibrium phase transition in $\ell^{2+\alpha}$, for any $\alpha>0$, will have sharp dynamic phase transition. In particular, \citet{berger/hoffman/sidoravicius/2005} introduced for any positive $\alpha$ a class of  kernels in $\ell^{2+\alpha}$ that exhibits equilibrium phase transition (see also \citet{dias/friedli/2015}). This implies that the same class of models exhibits dynamic phase transition for kernels in $\ell^{2+\alpha}$ and dynamic uniqueness for $\ell^2$. 

The coincidence between equilibrium and dynamic uniquenesses/phase transitions does not hold in general.  In what follows, we will exhibit two models studied in the literature that can exhibit equilibrium uniqueness and dynamic phase transition at the same time. 

\cite{friedli/2010} studied the following generalization of the model introduced by \cite{bramson/kalikow/1993}. Let $S = \{-1,+1\}$, $\vareps \in (0,1/2)$, and $(m_j)_{j \geq 1}$ be an increasing sequence of positive odd numbers. We consider a non-decreasing function $\psi:[-1, 1] \to [\varepsilon, 1-\varepsilon]$ which satisfies $\psi(r)+\psi(-r) = 1$. We call the sequence $(m_j)_{j \geq 1}$ \emph{lacunary} if for some $0\leq r_0 < 1$ such that $\psi(r_0) >  \psi(-r_0)$ we have $m_{j+1} \geq \frac{4}{1-r_0}m_{j}$ for all $j \geq 1$. Let $x \in \mathcal{X}$, for $j\geq 1$ we denote by $Q_{j}$ the function
\begin{equation*} 
\label{eq:BKpart}
Q_{j}(a\x) = \psi \bigg(\frac{a}{m_j}\sum_{l=1}^{m_j}x_{-l} \bigg).
\end{equation*}
Let $(\lambda_j)_{j \geq 1}$ be a sequence of positive numbers such that $\sum_{j=1}^\infty \lambda_j = 1$. Given $(m_j)_{j \geq 1}$ and $(\lambda_j)_{j \geq 1}$, the BKF-model associated to the parameters  $(m_j)_{j \geq 1}$ and $(\lambda_j)_{j \geq 1}$ is given by the kernel $g$ such that, for all $x \in \mathcal{X}$,
\begin{equation}
\label{eq:BKdef}
g(a\x) = \sum_{j = 1}^\infty \lambda_j Q_{j}(a\x).
\end{equation}
The original model introduced by \citet{bramson/kalikow/1993} is obtained by choosing $\psi(r) = 1-\varepsilon$ if $r \geq 0$ and $\psi(r) = \varepsilon$ otherwise. \citet{bramson/kalikow/1993} showed that it is possible to choose $(m_j)_{j \geq 1}$ and $(\lambda_j)_{j \geq 1}$ such that the corresponding model exhibits phase transition. Some progress has been made to obtain sufficient conditions for equilibrium phase transition in this model \citep{friedli/2010, gallo2014attractive, gallesco/gallo/takahashi/2014}, but a sharp condition on the kernel to guarantee equilibrium phase transition remains elusive. We note that in all known conditions of equilibrium phase transition of the BK model, the sequence $(m_j)_{j \geq 1}$ is lacunary \citep{friedli/2010, gallesco/gallo/takahashi/2014}.

 The binary autoregressive models constitute a different class of models that can also exhibit equilibrium phase transition. These models are defined through the following parameters: a continuous and increasing function $\phi:\mathbb{R}\rightarrow]0,1[$ such that $\phi(r)+\phi(-r) = 1$,  an absolute summable sequence of real numbers $(\beta_{n})_{n\geq1}$, and a real parameter $\delta$. The kernel $g$ of a  binary autoregressive model on the alphabet $\{-1,+1\}$ is given by 
\[
g(a\x):=\phi\left(a\sum_{n\geq1}\beta_{n}x_{-n}+a\delta\right).
\]
If $\phi$ is \emph{bi-Lipschitz}, we have that $1/\gamma\sum_{n>k}\beta_{n}\leq \textrm{var}_{k}(g) \leq \gamma\sum_{n>k}\beta_{n}$ for some positive constant $\gamma$. 
An important example of  binary autoregressive model is when $\phi(r)=e^{-r}(e^{-r}+e^{r})^{-1}$. The resulting kernel is called \emph{logit model} in the statistics literature, and \emph{one-sided $1$-dimensional long-range Ising model} in statistical physics literature. When  $(\beta_{n})_{n\geq1}$ and $\delta$ are non-negative real numbers, we say that the binary autoregressive model is attractive.  \citet{hulse/2006} showed that a class of attractive logit models exhibits equilibrium phase transition. 

We prove the following result.

\begin{coro}
\label{theo1}
A BKF model with lacunary $(m_j)_{j \geq 1}$ or an attractive binary autoregressive model with bi-Lipshitz $\phi$ exhibits dynamic uniqueness if and only if the corresponding  kernel satisfies the $\ell^2$ criterion.
\end{coro}
This result should be contrasted with the results for equilibrium phase transition for BKF models in the literature. For instance, the known explicit conditions for phase transition of Bramson-Kalikow model \citep{friedli/2010, gallesco/gallo/takahashi/2014} assume  kernels with continuity rates that are not in $\ell^p$ for any positive $p$. It is still an interesting open problem to decide whether there is any BKF or attractive binary autoregressive model with kernel in $\ell^p$, for some positive $p$, which exhibits equilibrium phase transition.
Let us define the \emph{oscillation} of order $k$ of a kernel $g$ as
\begin{equation*}
\textrm{osc}_{k}(g):=\sum_{a\in S}\sup_{\omega \in \X}\sup_{b, b' \in S}|g(a\omega^{-1}_{-k+1}b\omega^{-k-1}_{-\infty})-g(a\omega^{-1}_{-k+1}b'\omega^{-k-1}_{-\infty})|.
\end{equation*}
\citet{fernandez/maillard/2005} showed that if $\sum_{k\geq 1}\textrm{osc}_{k}(g) < 1$ (\emph{one-sided Dobrushin condition}), then there is a unique $g$-measure. For BKF and binary autoregressive models, it is not difficult to exhibit models that are not in $\ell^2$, but satisfy \citet{fernandez/maillard/2005} uniqueness criterion. For example, if $\psi(r) = 1/2+(1/2-\vareps)r$, the BKF model always satisfies \citet{fernandez/maillard/2005} uniqueness criterion although the kernels can have arbitrarly slow continuity rates. This shows that the equilibrium and dynamic uniquenesses/phase transitions are not equivalent in general. The implication is that there exists a kernel $g$ with a unique compatible measure $P$ and a non-stationary chain $P^{\x}$ (for some  $x \in \mathcal{X}$) such that $P$ and $P^{\x}$ differ on $\mathcal{T}^+$. This suggests the following physical interpretation: when dynamic and equilibrium uniqueness do not coincide, there is a ``hidden" phase transition (dynamic phase transition) that cannot be detected looking at events that depend on a finite number of coordinates, but can  only be detected looking at the asymptotic events.

To conclude the observations on Theorem \ref{theo:dynamicunique}, let us mention that, thanks to the ``$P^{\x}$-a.e.\ $\omega$'' in $(viii)$, it is not necessary that a kernel has a continuity rate in $\ell^2$ to exhibit dynamic uniqueness. As an example, let us construct a very simple kernel $g$ on $S=\{-1,+1\}$ having arbitrarily slow continuity rate and satisfying dynamic uniqueness. Let $(q_i)_{i\ge0}$ be a real sequence of $(0,1)$-valued numbers such that $q_i\searrow q_\infty>0$ and for any $\omega\in \mathcal X$ and $n\geq -1$, let $r(\omega_{-\infty}^n):=\sup\{k\leq n:\omega_{-k}=+1\}$ (by convention set $\sup\emptyset=-\infty$). The kernel for the example is defined by $g(+1\underline{\omega})=q_{r(\underline{\omega})}$. This is the kernel of a renewal sequence where the occurrence of symbol $1$ is a renewal event. A simple calculation shows that $\textrm{var}_k(g)=q_k-q_\infty$, therefore it is regular. Since $\inf_{\underline{x}}g(+1\underline{x})=q_\infty>0$, for any $x\in\mathcal{X}$ and for $P^{\x} \text{-a.e.}\; \omega$, there exists $N\in \Z_+$ such that $r(\omega_{-\infty}^N)\geq 0$.
Therefore, for all $\x,\y \in \mathcal{X}^-$,
\begin{align*}
\sum_{n=0}^\infty \sum_{a \in S}\left(g(a\omega^{n}_0\x)-g(a\omega^{n}_0\y)\right)^2 &=2\sum_{n=0}^\infty\left(g(+1\omega^{n}_0\x)-g(+1\omega^{n}_0\y)\right)^2\\
&=2\sum_{n=0}^\infty\left(q_{r(\omega^{n}_0\x)}-q_{r(\omega^{n}_0\y)}\right)^2\\
&\leq2( r(\omega_{-\infty}^N)+1)< \infty, \;\text{for}\; P^{\x} \text{-a.e.}\; \omega.
\end{align*}
Thus, $g$ exhibits dynamic uniqueness.
On the other hand since
$\textrm{var}_k(g)=q_k-q_\infty$ the continuity rate of $g$ is controlled by the rate of convergence of $q_k$ to its limit, which can be arbitrarily slow. \\

Theorem  \ref{theo:dynamicunique} shows that if a regular kernel satisfies dynamic uniqueness, the unique $g$-measure is mixing, or in the language of ergodic theory, it is a K-automorphism. But given that the dynamic uniqueness guarantees that all chains starting with different pasts are also mixing, it is natural to ask whether dynamic uniqueness implies some stronger form of mixing for the unique $g$-measure. We will show that this is indeed the case. Before stating the result we need the following definition. 

A stationary process $(X_j)_{j\in \Z}$ on $(\X, \F, \mu)$ or equivalently the measure $\mu$ is called \emph{weak Bernoulli} or, equivalently, \emph{absolute regular} or \emph{$\beta$-mixing} \citep{berbee/1986}, if 
\begin{equation}\lim_{n \rightarrow \infty} d_{TV}(\mu|_{\F^-}\otimes\mu|_{\F_{[n, \infty)}}, \mu|_{\F^-\otimes\F_{[n, \infty)}}) = 0.
\end{equation}
In the above limit, $\F^-\otimes\F_{[n, \infty)}$ is the sub $\sigma$-algebra of $\F$ generated by the rectangles $A\times B$ where $A\in \F^-$ and $B\in \F_{[n, \infty)}$ and $\mu|_{\F^-}\otimes\mu|_{\F_{[n, \infty)}}$ is the product measure between $\mu|_{\F^-}$ and $\mu|_{\F_{[n, \infty)}}$.
Weak Bernoullicity implies a strong form of mixing. To see this, let $\bigcap_{n \geq 1} \Big(\F_{(-\infty, -n]} \vee  \F_{[n, \infty)} \Big)$ be the two-sided $\sigma$-algebra. It is well known that if a stationary process is weak Bernoulli, then its past, future, and two-sided tail $\sigma$-algebras are trivial \citep{bradley/2005}.   

We prove the following result. 
\begin{theo} \label{theo:weakBernoulli}
 If a kernel $g$ is regular and exhibits dynamic uniqueness then the unique $g$-measure is weak Bernoulli. 
\end{theo}

An immediate consequence of the above theorem is that if a regular kernel $g$ is in $\ell^2$ then the unique $g$-measure is weak Bernoulli. We note that \citet{johansson/oberg/pollicott/2012} proved that if the kernel is in $\ell^2$, then the unique $g$-measure is very weak Bernoulli, which is weaker than being weak Bernoulli. 

  To prove Theorem \ref{theo:weakBernoulli} we will characterize all regular $g$-measures that are weak Bernoulli. The result has interest in its own and can be seen as a natural generalization of Theorem  \ref{theo:dynamicunique}. To emphasize this point, let us introduce the following definition, which is a generalization of the definition of dynamic uniqueness when there is a stationary measure compatible with the kernel.
   
  \begin{defi}
Given a regular kernel $g$, let $P$ be a $g$-measure. We say that the kernel $g$ exhibits \emph{$P$-dynamic uniqueness}~if for $P|_{\F^-} \otimes P|_{\F^-}$-a.e.~pair $(\x,\y)$ we have $P^{\x}|_{\mathcal{T}^+} = P^{\y}|_{\mathcal{T}^+}$, otherwise we say that there is \emph{$P$-dynamic phase transition}.
\end{defi}
Observe that in the above definition, the $g$-measure $P$ does not need to be unique. 
To state the next result, we will need to introduce one more definition: a $g$-measure is \emph{extremal} if  it cannot be written as a convex combination of other $g$-measures. The following theorem characterizes $g$-measures that are weak Bernoulli.

\begin{theo} \label{theo:betamixing}
Given a regular kernel $g$, let $P$ be an extremal $g$-measure. The following statements are equivalent:
\begin{enumerate}[(i)]
\item $P$ is weak Bernoulli (absolute regular, $\beta$-mixing). \label{00}

 \item $g$ exhibits $P$-dynamic uniqueness. \label{01}
 
 \item For $P|_{\F^-} \otimes P|_{\F^-}$-a.e.~$(\x,\y) $, \label{02}
$
 \lim_{n \rightarrow \infty} \sup_{B \in \F_{[n, \infty)}} |P^{\x}[B] - P^{\y}[B]| = 0.
$

\item For $P|_{\F^-} \otimes P|_{\F^-}$-a.e.~$(\x,\y) $ there exist a coupling $\tilde{P}$ between $P^{\x}$ and $P^{\y}$ and a coupling time $\Theta$ such that $\tilde{P}[\Theta < \infty] = 1$. \label{03}

 \item For $P|_{\F^-}$-a.e.~$\x$, $P^{\x}[A] = P^{\x}[A]^2 = P[A]$ for all $A \in \mathcal{T}^+$. \label{04}
 
\item For $P|_{\F^-} \otimes P|_{\F^-}$-a.e.~$(\x,\y)$, $P^{\x}|_{\F^+} \ll P^{\y}|_{\F^+}$. \label{05}

\item $g$ satisfies the $P$-weak-$\ell^2$ criterion.

\end{enumerate}
\end{theo}
For the above theorem, assuming that the $g$-measure is extremal is not a loss of generality because weak Bernoulli measures are trivial on the past tail $\sigma$-algebra and a $g$-measure is extremal if and only if  it is trivial on the past tail $\sigma$-algebra \citep{fernandez/maillard/2005, walters/2000}.

Besides proving Theorem \ref{theo:weakBernoulli}, the above theorem is useful for deriving some other sufficient conditions for weak Bernoullicity of $g$-measures. 

In the following corollary $L^\infty(\mathbf{T})$ stands for the usual Lebesgue space of complex valued functions on the unit circle $\mathbf{T}$. Given a probability measure $\mu$, we also use the notation $E_{\mu}$ for the expectation associated to $\mu$.
\begin{coro} \label{cor:crosscorrelation}
 Assume that $g$ is a binary autoregressive model with bi-Lipshitz $\phi$. Let $P$ be an extremal $g$-measure, $(\xi_j)_{j \in \Z}$ be the associated process, and $\gamma_j = E_P[\xi_0\xi_j]-E_P[\xi_0]E_P[\xi_j]$. If $(\gamma_j)_{j\in \Z}$ are the Fourier coefficients for some function $f \in L^\infty(\mathbf{T})$, \textit{i.e.},
 \begin{equation*}
 \gamma_j = \frac{1}{2\pi} \int^{2\pi}_{0}f(e^{i\theta})e^{-ij\theta}d\theta, \;\; j\in\Z, i = \sqrt{-1},
\end{equation*}
then $P$ is weak Bernoulli. In particular, if $\sum_{j \geq 0} |\gamma_j| < \infty$ then $P$ is weak Bernoulli.
\end{coro}
Observe that in the above corollary, we do not assume that the binary autoregressive kernel is attractive.

We can now give an example of a kernel $g$ that has dynamic phase transition but exhibits $P$-dynamic uniqueness. 
\begin{coro} \label{coro:betaexample}
 Let $g$ be a long range 1-dimensional Ising model.  If the parameters are such that $\delta \geq 0$ and  for $j \geq 1$, $\beta_j = c/j^{1+\epsilon}$, for a positive $c$ and some $\epsilon \in (0,1/2)$ with $\sum_{j = 1}^\infty \beta_j < 1$, then the $g$-measure is unique and weak Bernoulli (absolute regular, $\beta$-mixing) although $g$ exhibits dynamic phase transition.
\end{coro}

The above result shows that when there is a $g$-measure $P$, the definition of $P$-dynamic uniqueness is a proper generalization of the definition of dynamic uniqueness. \\

Let us summarize the results of this paper. Recall that all these results are about regular kernels. 
\begin{itemize}
\item Dynamic uniqueness $\Leftrightarrow$ convergence in total variation $\Leftrightarrow$ weak-$\ell^2$ criterion. 
\item weak-$\ell^2$ criterion $\Leftarrow \ell^2$ criterion.  
\item weak-$\ell^2$ criterion $\nRightarrow \ell^2$ criterion. 
\item Dynamic uniqueness $\Rightarrow$ equilibrium uniqueness.
\item Dynamic uniqueness $\nLeftarrow$ equilibrium uniqueness. 
\item For the Bramson-Kalikow-Friedli and autoregressive models, 
\[\textrm{Dynamic uniqueness}\,\,\,\Leftrightarrow\,\,\textrm{$\ell^2$ criterion}.
\]
\item Weak Bernoully ($\beta$-mixing) $\Leftrightarrow$ $P$-weak-$\ell^2$ criterion.
\item Dynamic uniqueness $\Rightarrow$ Weak Bernoully ($\beta$-mixing).
\item Dynamic uniqueness $\nLeftarrow$ Weak Bernoully ($\beta$-mixing).
\end{itemize}


\section{Proof of Theorem \ref{theo:dynamicunique}}

\subsection*{Proof of  equivalences $(i)$, $(ii)$, and $(iii)$ in Theorem \ref{theo:dynamicunique}}
These results are straightforward consequences of some well known results for stochastic processes.  
We state the general theorem from \citet{thorisson/2000} for the convenience of the reader.
\begin{theo*}[Thorisson (2000), Theorem 9.4, chapter 4] \label{theo:thorisson}
 Let $Y=(Y_j)_{j \geq 0}$ and $Z=(Z_j)_{j \geq 0}$ be canonical processes on probability spaces $(S^{\N}, \F^+, \mu)$ and $(S^{\N}, \F^+, \nu)$, respectively. The following are equivalent:
 
 \begin{enumerate}[(a)]
 
 \item For all tail events $A \in \mathcal{T}^+$ we have $\mu[A] = \nu[A]$.

\item $\lim_{n \rightarrow \infty} \sup_{B \in \F_{[n, \infty)}} |\mu[B] - \nu[B]| = 0.$

\item There exists a coupling $\tilde{P}$ between $Y$ and $Z$ such that the coupling time $\Theta$ is finite $\tilde{P}$-a.s.
 
 \end{enumerate}
 
\end{theo*}
Taking $\mu = P^{\x}|_{\F_{[0, \infty)}}$ and $\nu = P^{\y}|_{\F_{[0, \infty)}}$, the equivalences $(a)$, $(b)$, and $(c)$ in the above theorem implies, respectively, the equivalences between $(i), (ii)$, and $(iii)$ in Theorem \ref{theo:dynamicunique}. 
\qed

\subsection*{Proof of equivalences $(i)$ and $(iv)$  in Theorem \ref{theo:dynamicunique}}
Equivalence between  $(i)$ and $(iv)$ is a straightforward consequence of the representation, for all $\lambda \in \mathcal{M}$ and any $B \in \F$,
\begin{equation*}
 P^\lambda[B] = \int_{\X}P^{\x}[B]\lambda[dx].
\end{equation*}
\qed

\subsection*{Proof of equivalences $(i)$ and $(v)$  in Theorem \ref{theo:dynamicunique}}
Let us prove a lemma suggested by a remark in \citet{olshen/1971}. In what follows, equality $\mu$-a.s.~of $\sigma$-algebras $\mathcal{T}^-$ and $\mathcal{T}^+$ means that for all $A\in \mathcal{T}^-$ there exists $B\in \mathcal{T}^+$ such that $\mu[A\Delta B]=0$ and vice-versa. We remind the reader that we always consider a finite alphabet $S$ and $\X = S^\Z$.

\begin{lemma} \label{lemma:futurepast}
 Let $\mu$ be a stationary probability measure on $(\X, \F)$. We have that $\mathcal{T}^- = \mathcal{T}^+$, $\mu$-a.s. 
\end{lemma}
\noindent

This result seems to be known in the dynamical systems literature (we have been told that it is due to Pinsker), but we were not able to find a reference explicitly proving it, so for the sake of completeness, we give a proof of this fact. 

\noindent
\begin{proof}
From now on and until the end of this proof, we consider that $\F$ is complete (that is, it contains the class of $\mu$-null sets). Furthermore, we assume that all the $\sigma$-algebras considered in this proof also contain the class of $\mu$-null sets.
Lemma \ref{lemma:futurepast} is a consequence of Theorem 2 in \citet{rohlin/sinai/1961}. We will now restate their theorem using our notation. Let $h(T)$ be the entropy of $T$ and $H(\mathcal{A}\mid T^{-1}\mathcal{A})$ be the conditional entropy of the measurable partition corresponding to $\mathcal{A}$ given $T^{-1}\mathcal{A}$. We also denote by $\pi(T)$ the Pinsker $\sigma$-algebra defined by $\pi(T):= \{A \in \F: h(\sig(A), T) = 0\}$, where $h(\sig(A), T)$ is the entropy of $T$ given the measurable partition corresponding to $\sig(A)$ ($\sig(A)$ denotes the $\sigma$-algebra generated by~$A$). For a definition of these notions, see for example \cite{martin/england/1981}.

\begin{theo*}[\citet{rohlin/sinai/1961}, Theorem 2]
Let $(\X, \F, \mu)$ be a Lebesgue probability space and $T$ be a measure preserving automorphism. If there exists a sub $\sigma$-algebra $\mathcal{A} \subset \F$ such that 
\begin{enumerate}[(a)]
 \item $T^{-1}\mathcal{A} \subset \mathcal{A}$,
 
 \item $\bigvee_{n=0}^{\infty}T^n\mathcal{A} = \F $,
 
 \item $h(T) = H(\mathcal{A} \mid T^{-1}\mathcal{A})$,
 
 \item $h(T)  < \infty$,
\end{enumerate}
then $\bigcap_{n=0}^{\infty}T^{-n}\mathcal{A} = \pi(T)$.
 
\end{theo*}

 To prove Lemma \ref{lemma:futurepast}, we will first verify that all the conditions of the above theorem are satisfied taking  $\mathcal{A} = \F_{[0, \infty)}$ and $T$ the shift defined by $(Tx_j) = x_{j+1}$. We have $T\F_{[j, \infty)} = \F_{[j-1, \infty)}$ and $T^{-1}\F_{[j, \infty)} = \F_{[j+1, \infty)}$, therefore, $(a)$ and $(b)$ in the above theorem are satisfied. By Theorem 2.39 in \cite{martin/england/1981}, we have that $h(T) = h(\F_{\{0\}}, T)$. Besides, since $\F_{\{0\}}$ is finite (up to $\mu$-null sets), by Theorem 2.27 in \cite{martin/england/1981}, we have 
 $$h(\F_{\{0\}}, T)=H(\F_{\{0\}}\mid \F_{[1,\infty)}).$$
 Since $H(\F_{\{0\}}\mid \F_{[1,\infty)})\leq \log |S|$, $(d)$ is satisfied. Finally, we have that $H(\F_{[0,\infty)}\mid \F_{[1,\infty)})=H(\F_{\{0\}}\mid \F_{[1,\infty)})$, which implies $(c)$. Thus using the above theorem, we conclude that $\mathcal{T}^+ = \pi(T)$. Similarly, if we take $\tilde{T}$ as the inverse shift $(\tilde{T}x_j) = x_{j-1}$ and  $\mathcal{A} = \F_{(-\infty, 0]}$, we conclude that $\mathcal{T}^- = \pi(\tilde{T})$. Now, it remains to show that $\pi(T) = \pi(\tilde{T})$. For this, we note that, if $T$ is invertible, we have (see again Theorem 2.27 in \cite{martin/england/1981})
\begin{align*}
 h(\sig(A), T) &= H(\sig(A) | \bigvee_{j=1}^{\infty}T^{-j}\sig(A)) \\
 &= H(\sig(A) | \bigvee_{j=1}^{\infty}T^j\sig(A))\\
 &= h(\sig (A), \tilde{T}). 
\end{align*}
Therefore, $\mathcal{T}^+ = \mathcal{T}^-$, as we wanted to show.
\end{proof}

Now, let us prove the equivalences $(i)$ and $(v)$ in Theorem \ref{theo:dynamicunique}. To prove that $(v)$ implies $(i)$, let $\lambda, \lambda' \in \mathcal{M}$ and take $\lambda'' = (\lambda+\lambda')/2$. For any $A \in \mathcal{T}^+$, we have $P^{\lambda''}[A] = (P^{\lambda}[A]+P^{\lambda'}[A])/2$. Since $P^{\lambda''}$ is trivial on $\mathcal{T}^+$, then $P^{\lambda}[A] = P^{\lambda'}[A] = P^{\lambda''}[A] =  1$ or $0$, as we wanted to show. 

To show that $(i)$ implies $(v)$, we observe that because $(i)$ implies $(iv)$, the existence of a $P \in \mathcal{M}_0$ trivial on $\mathcal{T}^+$  implies that all elements of $\mathcal{M}_0$ are trivial. To prove that such $P$ always exists, we first use the fact that extremal $g$-measures are trivial on $\mathcal{T}^-$ \citep{fernandez/maillard/2005, walters/2000} and that there is  at least one extremal $g$-measure \citep{walters/2000}. Therefore, we conclude that the extremal $g$-measures are also trivial on  $\mathcal{T}^+$. Second, let $P$ be an extremal $g$-measure and $P_0 = P|_{\F_{(-\infty, -1]}}$. The stationary measure $P$ can be written as $P^{P_0}$ and therefore it is an element of $\mathcal{M}_0$, as we wanted to show.
\qed

\subsection*{Proof of equivalences $(v)$ and $(vi)$  in Theorem \ref{theo:dynamicunique}}
The equivalence between $(v)$ and $(vi)$ is standard in the literature (see for example Theorem 2 in \citet{blackwell/freedman/1964}). We exhibit the argument for convenience of the reader.

If $(vi)$ holds, we have for any $A \in \mathcal{T}^+$, $P^\lambda[A] = P^\lambda[A]^2$. Hence, $P^\lambda[A] = 0$ or $1$. To prove that $(v)$ implies $(vi)$, fix $B \in \F$. For any $A \in \F_{[n, \infty)}$, we have
\begin{align*}
 \big|P^\lambda[A \cap B] -P^\lambda[A] P^\lambda[B]\big| &= \big|\int_{A}  (P^\lambda[B | \F_{[n, \infty)}] - P^\lambda[B])dP^\lambda\big| \\
 & \leq \int_{\X}  \big|P^\lambda[B | \F_{[n, \infty)}] - P^\lambda[B] \big|dP^\lambda.
\end{align*}
By backward martingale convergence theorem and triviality of $P^\lambda$ in $\mathcal{T}^+$, the right hand side of the above inequality converges to $0$. Therefore, we have
\begin{equation*}
\lim_{n \to \infty} \sup_{A\in \F_{[n, \infty)}} \big|P^\lambda[A \cap B] -P^\lambda[A] P^\lambda[B]\big| = 0,
\end{equation*}
as we wanted to show.
\qed

\subsection*{Proof of equivalences $(i)$ and $(vii)$ in Theorem \ref{theo:dynamicunique}}
 \cite{engelbert/shiryaev/1980} proved the following dichotomy result.
\begin{theo*}[Engelbert \& Shiryaev (1980), Theorem 5] \label{theo:dicotomy}
Let $Y=(Y_j)_{j \geq 0}$ and $Z=(Z_j)_{j \geq 0}$ be canonical processes on probability spaces $(S^{\N}, \F^+, \mu)$ and $(S^{\N}, \F^+, \nu)$, respectively. If for all $n \geq 0$, we have $\mu |_{\F_{[0,n]}} \ll \nu |_{\F_{[0,n]}}$ then
\begin{enumerate}[(a)]
\item $\mu \ll \nu$ if and only if $$\mu \left[\limsup_{n \rightarrow \infty}\frac{d\mu|_{\F_{[0,n]}}}{d\nu|_{\F_{[0,n]}}} < \infty \right] = 1.$$

\item $\mu \perp \nu$ if and only if $$\mu \left[\limsup_{n \rightarrow \infty}\frac{d\mu|_{\F_{[0,n]}}}{d\nu|_{\F_{[0,n]}}} < \infty \right] = 0.$$

\end{enumerate}

\end{theo*}

To prove that $(i)$ implies $(vii)$, let $\mu = P^{\x}|_{\F^+}$ and $\nu = P^{\y}|_{\F^+}$ in the above theorem.  The regularity of kernel $g$ guarantees $P^{\x}|_{\F_{[0,n]}} \ll P^{\y}|_{\F_{[0,n]}}$ for all $n \geq 0$.  We observe that the event
\begin{equation*}
 D = \left \{ \limsup_{n \rightarrow \infty}\frac{dP^{\x}|_{\F_{[0,n]}}}{dP^{\y}|_{\F_{[0,n]}}} < \infty  \right\}
\end{equation*}
is in $\mathcal{T}^+$. 
We have that $(i)$ implies $(v)$ in Theorem \ref{theo:dynamicunique}. Therefore, for each $\x, \y \in \X^-$, we have $P^{\x}[D] = P^{\y}[D] = 0$ or $1$. By Fatou's  lemma, $E_{P^{\y}}\Big[ \limsup_{n \rightarrow \infty}\frac{dP^{\x}|_{\F_{[0,n]}}}{dP^{\y}|_{\F_{[0,n]}}} \Big] \leq 1$ and therefore $P^{\y}[D] = 1$. Hence, we have $P^{\x}[D] = 1$.  As a consequence, the above theorem by \citet{engelbert/shiryaev/1980} implies that $P^{\x}|_{\F^+} \ll P^{\y}|_{\F^+}$ as we wanted to show.

To prove that $(vii)$ implies $(i)$, we use again the fact that an extremal $g$-measure $P$ is trivial on $\mathcal{T}^-$ and that there is  at least one extremal $g$-measure. Because of Lemma  \ref{lemma:futurepast}, this implies that $P$ is trivial on $\mathcal{T}^+$. Using the property of regular conditional probability, we have that for any $A \in \mathcal{T}^+$, $P^{\x}[A] = P[A | \F^-](x)=0$ or $1$ for $P$-a.e.~$x$.  In particular, for any $A$ such that $P[A] = 0$, there is a past $\underline{\omega}_A \in \X^-$ such that $P^{\underline{\omega}_A}[A] = 0$.  If $(vii)$ holds,  for all $\y \in \X^-$ and $A \in \mathcal{T}^+$, we have that $P^{\y}[A] =  P^{\underline{\omega}_A }[A]= 0$. This implies that, for all $\x, \y \in \X^-$, $P^{\x}|_{\mathcal{T}^+} = P^{\y}|_{\mathcal{T}^+}$.
\qed

\subsection*{Proof of equivalence $(vii)$ and $(viii)$ in Theorem \ref{theo:dynamicunique}}
To prove the equivalence between  $(vii)$ and $(viii)$, we use the results from the theory of predictable absolute continuity and singularity (ACS) criteria developed by Shiryaev and co-authors. The idea of using the predictable ACS criteria for $g$-chains was initiated in \citet{johansson/oberg/pollicot/2007}. In what comes next, we closely follow the exposition in \citet{johansson/oberg/pollicot/2007}.

Let $X=(X_j)_{j \in \Z}$ and $Y=(Y_j)_{j \in \Z}$ be canonical processes on $(\X, \F, \mu)$ and $(\X, \F, \nu)$, respectively. For $\omega \in \X$, let $Z_n(\omega) := \frac{d\mu|_{\F_{[0,n]}}}{d\nu|_{\F_{[0,n]}}}(\omega)$ and $\alpha_n(\omega) := Z_n(\omega)/Z_{n-1}(\omega)$. Also, define
\begin{equation*}
 d_n(\omega) := E_{\nu}\big[(1-\sqrt{\alpha_n})^2 | \F_{[0, n-1]}\big](\omega). 
\end{equation*}

The predictable ACS criteria is given by
\begin{theo*}[see \cite{jacod/shiryaev/2002}, Theorem 2.36, p.253]
If for all $n \geq 0$ we have $\mu |_{\F_{[0,n]}} \ll \nu |_{\F_{[0,n]}}$, then $\mu|_{\F^+} \ll \nu|_{\F^+}$ if and only if $\sum_{n=1}^\infty d_n(\omega) < \infty$, $\mu$-a.s.
\end{theo*}

Let us rewrite $d_n(x)$ in a more explicit form. We have that
\begin{equation*}
 \frac{d \mu|_{\F_{[0,n]}}}{d \mu|_{\F_{[0,n-1]}}}(\omega) = \mu(X_n = \omega_n | X^{n-1}_0=\omega^{n-1}_0)=: \mu(\omega_n | \omega^{n-1}_0).
\end{equation*}

Hence, 
\begin{align*}
 E_{\nu}\big[(1-\sqrt{\alpha_n})^2 | \F_{[0, n-1]}\big](\omega) &= \sum_{\omega_n \in S} \left(1-\sqrt{\frac{\mu(\omega_n | \omega^{n-1}_0)}{\nu(\omega_n | \omega^{n-1}_0)}}\right)^2 \nu(\omega_n | \omega^{n-1}_0)\\
 &= \sum_{\omega_n \in S} \left(\sqrt{\nu(\omega_n |\omega^{n-1}_0)}-\sqrt{\mu(\omega_n | \omega^{n-1}_0)}\right) ^2.
\end{align*}

To prove the equivalence $(vii)$ and $(viii)$ in Theorem \ref{theo:dynamicunique}, we take $\mu = P^{\x}$, $\nu = P^{\y}$ and apply the above theorem. Observe that by the strong non-nullness assumption (\ref{strongnonnull}), we have, for all $n \geq 0$, that $P^{\x} |_{\F_{[0,n]}} \ll P^{\y} |_{\F_{[0,n]}}$. In our case, we have
\begin{align*}
 d_n(\omega) &= \sum_{\omega_n \in S} \left(\sqrt{P^{\y}(\omega_n | \omega^{n-1}_0)}-\sqrt{P^{\x}(\omega_n | \omega^{n-1}_0)}\right)^2 \\
 &= \sum_{a \in S} \left(\sqrt{g(a \omega^{n-1}_0 \y)}-\sqrt{g(a \omega^{n-1}_0 \x)}\right)^2.
\end{align*}

To conclude, we only need to show that $\sum_{a \in S} (g(a \omega^{n-1}_0 \y)^{1/2}-g(a \omega^{n-1}_0 \x)^{1/2})^2$ has the same order as $\sum_{a \in S} (g(a \omega^{n-1}_0 \y)-g(a \omega^{n-1}_0 \x))^2$. We have that
\begin{align*}
 &\sum_{a \in S} (g(a \omega^{n-1}_0 \y)-g(a \omega^{n-1}_0 \x))^2  \\
 &= \sum_{a \in S} (g(a \omega^{n-1}_0 \y)^{1/2}-g(a \omega^{n-1}_0 \x)^{1/2})^2(g(a \omega^{n-1}_0 \y)^{1/2}+g(a \omega^{n-1}_0 \x)^{1/2})^2.
\end{align*}
Therefore, taking 
$$\gamma = \inf_{a \in S}\inf_{\underline{\omega} \in \X^-} g(a \underline{\omega}),$$
which belongs to the interval $(0,1)$ by (\ref{strongnonnull}),  we have
\begin{equation*}
 4\gamma d_n(\omega) \leq \sum_{a \in S} (g(a \omega^{n-1}_0 \y)-g(a \omega^{n-1}_0 \x))^2 \leq 4(1-\gamma)d_n(\omega).
\end{equation*}

From the above inequality and the ACS criteria, we conclude that $P^{\x}|_{\F^+} \ll P^{\y}|_{\F^+}$ if and only if 
\begin{equation*}
 \sum_{n=0}^\infty\sum_{a \in S} (g(a \omega^{n}_0 \y)-g(a \omega^{n}_0 \x))^2 < \infty
\end{equation*}
for $P^{\x}$-a.e. $\omega$.
\qed

\section{Proof of Corollaries \ref{coro:dynamicphase},  \ref{theo:dynamicuniquel2}, and \ref{theo1}}

\subsection*{Proof of Corollary \ref{coro:dynamicphase}}
If there is equilibrium phase transition, there exist two stationary probability measures $P$ and $P'$ compatible with $g$, such that $P[C] \neq P'[C]$ for some cylinder $C$. Because $P$ and $P'$ are stationary 
\begin{equation}
\label{Contra}
|P[T^{n}C] - P'[T^{n}C]| = |P[C] - P'[C]| > 0
\end{equation} for any $n \in \Z$.
Now, if there is dynamic uniqueness, by Theorem~\ref{theo:dynamicunique} we have that $(\ref{4})$ holds for $P$ and $P'$.  Using again the theorem by \citet{thorisson/2000}, we have
\begin{equation*}
 \lim_{n \rightarrow \infty} \sup_{B \in \F_{[n, \infty)}} |P[B] -P'[B]| = 0,
\end{equation*}
which contradicts (\ref{Contra}). Thus, $g$ exhibits dynamic phase transition.
\qed

\subsection*{Proof of Corollary \ref{theo:dynamicuniquel2}}
We have by definition of variation rate that for all $a \in S$, $\omega \in \X$ and $\x, \y \in \X^-$ 
\begin{equation*}
 |g(a \omega^{n}_0 \y)-g(a \omega^{n}_0 \x)| \leq \textrm{var}_{n}(g),
\end{equation*}
therefore, 
\begin{equation*}
 \sum_{n=0}^\infty\sum_{a \in S} (g(a \omega^{n}_0 \y)-g(a \omega^{n}_0 \x))^2 \leq |S| \sum_{n=0}^\infty \textrm{var}_n(g)^2.
\end{equation*}
Using Theorem \ref{theo:dynamicunique} (\textit{\ref{8}}), we conclude the result. 
\qed

\subsection*{Proof of Corollary \ref{theo1}}
Corollary \ref{theo:dynamicuniquel2} shows that if the kernel is in $\ell^2$  then we have dynamic uniqueness. Hence, from Theorem \ref{theo:dynamicunique} (\textit{\ref{8}}),  we only need to show that if the BKF and attractive binary autoregressive kernels are not in $\ell^2$, then there exist $\x, \y \in \X^-$  such that 
\begin{equation} \label{eq:diverge}
 \sum_{n=0}^\infty\sum_{a \in S} (g(a \omega^{n}_0 \y)-g(a \omega^{n}_0 \x))^2 = \infty
\end{equation}
for all $\omega\in \X$.

For both models, we choose $\x = \underline{1}$, where $\underline{1}_j = 1$ for all $j \leq -1$. Analogously, we define $\underline{-1}_j = -1$ for all $j \leq -1$ and we take $\y = \underline{-1}$. For both BKF and binary autoregressive models we have, for all $\omega \in \X$,
\begin{equation*}
 \sum_{n=0}^\infty\sum_{a \in S} (g(a \omega^{n}_0 \underline{-1})-g(a \omega^{n}_0 \underline{1}))^2  \geq \sum_{n=0}^\infty\inf_{a \in S} \inf_{\omega^n_0 \in S^{n}}(g(a \omega^{n}_0 \underline{-1})-g(a \omega^{n}_0 \underline{1}))^2.
\end{equation*}

A straightforward calculation shows that for the BKF model with lacunary $(m_j)_{j \geq 1}$, we have for any $m_n < j \leq \frac{(1-r_0)}{2}m_{n+1}$
\begin{equation*}
 \inf_{a \in S} \inf_{\omega^j_0 \in S^{j}}|g(a \omega^{j}_0 \underline{-1})-g(a \omega^{j}_0 \underline{1})| \geq (\psi(r_0)-\psi(-r_0))\sum_{k \geq n+1}\lambda_k.
\end{equation*}
Therefore, 
\begin{equation} 
\label{eq:lower}
\sum_{n=0}^\infty\inf_{a \in S} \inf_{\omega^n_0 \in S^{n}}(g(a \omega^{n}_0 \underline{-1})-g(a \omega^{n}_0 \underline{1}))^2 \geq \frac{1-r_0}{4}(\psi(r_0)-\psi(-r_0))^2  \sum_{n=2}^\infty m_{n}\big(\sum_{k \geq n}\lambda_k\big)^2.
\end{equation}
We also have
\begin{equation}
\label{eq:upper}
\sum_{n=1}^\infty \textrm{var}_n(g)^2 \leq (1-\varepsilon)^2 \sum_{n=1}^\infty m_{n} \big(\sum_{k \geq n}\lambda_k\big)^2.
\end{equation}
Therefore, from \eqref{eq:upper} we conclude that if $\sum_{n=1}^\infty \textrm{var}_n(g)^2$ diverges then $\sum_{n=1}^\infty m_{n} \big(\sum_{k \geq n}\lambda_k\big)^2$ diverges. Furthermore, from \eqref{eq:lower} if $\sum_{n=1}^\infty m_{n} \big(\sum_{k \geq n}\lambda_k\big)^2$ diverges then we obtain \eqref{eq:diverge}, as we wanted to show.

For the binary autoregressive model, using the fact that $\phi$ is bi-Lipschitz, we have
\begin{equation} \label{eq:lowerBin}
\sum_{n=0}^\infty\inf_{a \in S} \inf_{\omega^n_0 \in S^{n}}(g(a \omega^{n}_0 \underline{-1})-g(a \omega^{n}_0 \underline{1}))^2 \geq 1/\gamma^2 \sum_{n=1}^\infty \big(\sum_{k>n}\beta_{k}\big)^2.
\end{equation}
Also,
\begin{equation} \label{eq:upperBin}
\sum_{n=1}^\infty \textrm{var}_n(g)^2 \leq \gamma^2 \sum_{n=1}^\infty \big(\sum_{k>n}\beta_{k}\big)^2.
\end{equation}
Therefore, from  \eqref{eq:lowerBin}, \eqref{eq:upperBin}, and Theorem \ref{theo:dynamicunique} (\textit{\ref{8}}),  we conclude that the binary autoregressive process is in $\ell^2$ if and only if it has dynamic uniqueness as we wanted to prove.
\qed

\section*{Proof of  Theorems \ref{theo:weakBernoulli} and \ref{theo:betamixing}}
Theorem  \ref{theo:weakBernoulli} is an immediate consequence of Theorem \ref{theo:betamixing}.

\subsection*{Theorem \ref{theo:betamixing}}

The equivalence between $(i)$ and $(iii)$ is a consequence of the well known identity
\begin{equation} \label{eq:equivbeta}
d_{TV}(P|_{\F^-}\otimes P|_{\F_{[n, \infty)}}, P|_{\F^-\otimes\F_{[n, \infty)}}) = 2E_P\left [\sup_{B \in \F_{[n, \infty)}} |P[B|\F^-] - P[B]| \right]
\end{equation}
where $P[\cdot |\F^-]$ is any regular version of the conditional probability.
For a proof of the above equation, see \citet{volkonskii/rozanov/1961}, Lemma 4.1 (see also \citet{tong/vanrandel/2014}, Lemma 2.9). Observe that the application $\x\mapsto \sup_{B \in \F_{[n, \infty)}} |P[B|\F^-](\x) - P[B]|$ is measurable since the ``$\sup$" can be taken over the countable algebra of cylinder sets in~$\F_{[n, \infty)}$.

 Because $\sup_{B \in \F_{[n, \infty)}} |P[B|\F^-](\x) - P[B]|$ is pointwise decreasing in $n$, we have that the right-hand side of \eqref{eq:equivbeta} converges to zero if and only if $\sup_{B \in \F_{[n, \infty)}} |P[B|\F^-](\x) - P[B]|$ converges to zero for $P|_{\F^-}$-a.e $\x$. Therefore, $P$ is weak Bernoulli if and only if for $P|_{\F^-}$-a.e. $\x$
\begin{equation*}
\lim_{n\rightarrow \infty} \sup_{B \in \F_{[n, \infty)}} |P^{\x}[B] - P[B]| = 0.
\end{equation*}
Now, we want to show that the above equality is equivalent to $\lim_{n\rightarrow \infty} \sup_{B \in \F_{[n, \infty)}} |P^{\x}[B] - P^{\y}[B]| = 0$ for $P|_{\F^-} \otimes P|_{\F^-}$-a.e.~$(\x, \y)$. This is a straightforward consequence of the following bounds. 
\begin{equation*}
\sup_{B \in \F_{[n, \infty)}} |P^{\x}[B] - P^{\y}[B]| \leq \sup_{B \in \F_{[n, \infty)}} |P^{\x}[B] - P[B]| + \sup_{B \in \F_{[n, \infty)}} |P^{\y}[B] - P[B]|
\end{equation*}
and
\begin{equation}
\sup_{B \in \F_{[n, \infty)}} |P^{\x}[B] - P[B]| \leq \int_{\X} \sup_{B \in \F_{[n, \infty)}} |P^{\x}[B] - P^{\y}[B]|P(dy). \label{eq:bound1P}
\end{equation}
Therefore, we have shown the equivalence between $(i)$ and $(iii)$. The equivalence between $(ii)$ and $(iii)$, and $(ii)$ and $(iv)$ follows from Thorisson's theorem stated in the proof of Theorem \ref{theo:dynamicunique}. 

$(v)$ trivially implies $(ii)$. We already showed that $(ii)$ implies $(iii)$. Using inequality \eqref{eq:bound1P}, if $(iii)$ holds then for $P|_{\F^-}$-a.e. $\x$ we have $P|_{\mathcal{T}^{+}} = P^{\x}|_{\mathcal{T}^{+}}$. Because, $P$ is extremal, $P$ is trivial on $\mathcal{T}^-$, thus using Lemma \ref{lemma:futurepast}, $P$ is trivial on $\mathcal{T}^+$. Therefore, $P^{\x}$ will be trivial for $P|_{\F^-}$ a.e. $\x$, and we conclude that $(iii)$ implies $(v)$, as we wanted to show.

$(vi)$ and $(vii)$ follows from the proof of equivalence between $(vii)$ and $(viii)$ in Theorem \ref{theo:dynamicunique}. Therefore, it remains to show the equivalence between $(i)$ and $(vi)$. We already proved that $(i)$ implies $(v)$, and clearly $(v)$ implies $(vi)$. To prove that $(vi)$ implies $(i)$, we use Corollary 2.8 in \citet{tong/vanrandel/2014}. In our notation, it states that $P$ is weak Bernoulli if and only if for $P|_{\F^-} \otimes P|_{\F^-}$-a.e. $(\x,\y)$, there exist a $k \geq 0$ such that $P^{\x}$ and $P^{\y}$ are not mutually singular on $F_{[k,\infty)}$. Condition $(vi)$ implies that for $P|_{\F^-} \otimes P|_{\F^-}$-a.e. $(\x,\y)$, $P^{\x}$ and $P^{\y}$ are not mutually singular on $F_{[0,\infty)}$, hence $P$ is weak Bernoulli as we wanted to show.
 \qed

\section*{Proof of Corollary \ref{cor:crosscorrelation}}
Let $(\eta_j)_{j \in \Z}$ be an independent copy of $(\xi_j)_{j\in \Z}$. From Theorem \ref{theo:betamixing} $(vii)$ a sufficient condition for $P$ to be weak Bernoulli is that

\begin{equation*} 
E_{P\otimes P}\Big[\sum_{n=1}^\infty \big(\sum_{k>n}\beta_{k}(\xi_k-\eta_k)\big)^2\Big] < \infty.
\end{equation*}
By Tonelli's theorem, we can interchange the summation and integration on the left hand side of the above equation, and the above condition is equivalent to
\begin{equation*} 
\sum_{n=1}^\infty \Big(E_P\Big[\big(\sum_{k>n}\beta_{k}\xi_k\big)^2\Big]-E_P\Big[\sum_{k>n}\beta_{k}\xi_k\Big]^2\Big)< \infty.
\end{equation*}
Now, because $\sum_{j \geq 1}|\beta_j| < \infty$, using the dominated convergence theorem, we can write the above inequality as
\begin{equation*}
\sum_{n=1}^\infty \sum_{j > n}\sum_{k>n}\beta_j\beta_k \gamma_{j-k}< \infty.
\end{equation*}
where we recall that $\gamma_j$ are the correlation coefficients.
Let $\Gamma=(\Gamma_{jk})_{j,k\ge1}$ be the operator defined by $\Gamma_{jk} = \gamma_{j-k}$. This is a Toeplitz operator and it is a classical result that $\Gamma$ is bounded in $\ell^2(\Z_+)$ if and only if $(\gamma_j)_{j \in \Z}$ are the Fourier coefficients of a function $f \in L^\infty(\mathbf{T})$ (see for example Chapter 1 of \cite{bottcher/grudsky/1991}). 
If $\Gamma$ is bounded in $\ell^2(\Z_+)$, we have that
\begin{equation*}
\sum_{j > n}\sum_{k>n}\beta_j\beta_k \gamma_{j-k} \leq C\sum_{k>n}\beta_k^2
\end{equation*}
where $C$ is some positive constant. In particular, if $\sum_{j \geq 0} |\gamma_j| < \infty$, we can take $C = \sum_{j \geq 0} |\gamma_j|$.

To conclude, we only need to show that if $\Gamma$ is bounded in $\ell^2(\Z_+)$ and $\sum_{j \geq 1}|\beta_j| < \infty$, then $\sum_{n=1}^\infty \sum_{k>n}\beta_k^2< \infty$. Since $\sum_{j \geq 1}|\beta_j| < \infty$ the sequence $(n|\beta_n|)_{n\geq 1}$ is bounded by a constant $C'$. This implies that
\[
\sum_{n=1}^\infty \sum_{k>n}\beta_k^2=\sum_{n=1}^\infty n\beta_n^2=\sum_{n=1}^\infty n|\beta_n|.|\beta_n|\leq C'\sum_{n=1}^\infty|\beta_n|< \infty,
\]
as we wanted to show.
\qed

\section*{Proof of Corollary \ref{coro:betaexample}} 

For the proof, assume $g$ is the kernel satisfying the conditions of the Corollary \ref{coro:betaexample}. From  \eqref{eq:lowerBin}, we have that if $\beta_j = c/j^{1+\epsilon}$, with $\epsilon \in (0,1/2)$, then the kernel $g$ is not in $\ell^2$. Applying Corollary \ref{theo1}, we deduce that $g$ exhibits dynamic phase transition. From \citet{fernandez/maillard/2005}, a regular $g$-measure with $\sum_{k\geq 1}\textrm{osc}_{k}(g) < 1$ (one-sided Dobrushin condition) has exactly one compatible stationary measure, which is by definition extremal. Because $\sum_{k\geq 1}\textrm{osc}_{k}(g) \leq \sum_{j \geq 1}\beta_j < 1$, we conclude that there exists exactly one $g$-measure $P$.  Therefore, we only need to show that $g$ exhibits $P$-dynamic uniqueness. For this, we will first show that the two-point correlation function for $P$ has a summable decay. The required bound is a consequence of Corollary 5.21 in \citet{fernandez/maillard/2005}, where the authors obtain an upper bound for the correlation between bounded oscillation function when $g$ satisfies the one-sided Dobrushin condition. The upper bound in \citet{fernandez/maillard/2005} depends on the behavior of the matrix $D = \sum_{n = 1}^\infty \alpha^n$ where $\alpha=(\alpha_{i,j})_{(i,j)\in \Z^2}$ is defined by $\alpha_{ij}= \textrm{osc}_{i-j}(g)$ for $i>j$ and  $\alpha_{ij}= 0$ otherwise. For our proof, we only need a bound for two-point correlation, hence we state a specialized version of their result below. 

\begin{lemma}[\citet{fernandez/maillard/2005}, Corollary 5.21]
 Assume that $g$ is a regular kernel satisfying the one-sided Dobrushin condition. Let $(\xi_j)_{j \in \Z}$ be the stationary process compatible with $g$, then we have, for $m \leq l$, 
\begin{equation} 
\label{eq:twopointcorr}
\big|E_P[\xi_l \xi_m]-E_P[\xi_l]E_P[\xi_m]\big| \leq KA_{lm},
\end{equation}
where $K$ is a positive constant and $A_{lm} = D_{lm} + \sum_{n \leq m}D_{ln}D_{mn}$.
\end{lemma}

To obtain an upper bound for $A_{lm}$, we will use a result by \citet{jaffard/1990} (Proposition 3), who shows that if $M: \ell^2(\Z_+) \to \ell^2(\Z_+)$ is an invertible matrix with entries satisfying $|M_{ij}| \leq c_1(1+|i-j|)^{-(1+\eta)}$, for some positive constants $c_1$ and $\eta$, then its inverse $M^{-1}$ has entries satisfying $|M^{-1}_{ij}| \leq C_1(1+|i-j|)^{-(1+\eta)}$, for some positive constant $C_1$. First, observe that the matrix $A:=I-\alpha$ (where $I$ is the identity matrix) satisfies the conditions of the theorem (taking $\eta=\epsilon$) since $\alpha$ defines a bounded operator from $\ell^2(\Z_+) \to \ell^2(\Z_+)$ and $\|\alpha\|<1$ (here $\|\cdot\|$ stands for the operator norm). Then, by definition, $D = \sum_{n = 1}^\infty \alpha^n=A^{-1}-I$. Therefore, we conclude that $|D_{ij}| \leq C'(1+|i-j|)^{-(1+\epsilon)}$, for some positive constant $C'$. Using this upper bound on \eqref{eq:twopointcorr}, we obtain
\begin{align*}
 &\big|E_P[\xi_l \xi_m]-E_P[\xi_l]E_P[\xi_m]\big| \\
 &\leq K\big[C'(1+|l-m|)^{-(1+\epsilon)} + \sum_{n \leq m}C'(1+|l-n|)^{-(1+\epsilon)}C'(1+|m-n|)^{-(1+\epsilon)}\big]\\
 &\leq KC'(1+|l-m|)^{-(1+\epsilon)}+K{C'}^2(1+|l-m|)^{-(1+\epsilon)}\sum_{j \geq 0}(1+j)^{-(1+\epsilon)}\\
 &\leq \kappa (1+|l-m|)^{-(1+\epsilon)},
\end{align*}
where $\kappa := K(C'+ C'^2)\sum_{j \geq 0}(1+j)^{-(1+\epsilon)}$.

The above upper bound implies that $\sum_{j\geq 0}|\gamma_{j}|< \infty$, and using Corollary \ref{cor:crosscorrelation} we conclude that $P$ is weak Bernoulli  as we wanted to show.
\qed

\bigskip
\footnotesize
\noindent\textit{Acknowledgments.}
This work is supported by USP project ``Mathematics, computation, language and the brain'', FAPESP project ``NeuroMat'' (grant 2011/51350-6),  CNPq projects ``Stochastic Modeling of the Brain Activity'' (grant 480108/2012-9), ``Codifica\c c\~ao Finit\'aria e cadeias de longo alcance'' (grant 462064/2014-0), and ``Sistemas Estoc\'asticos: equil\'ibrio e n\~ao equil\'ibrio, limites em escala e percola\c c\~ao'' (grant 474233/2012-0). CG was supported by FAPESP (grant 2013/10101-9). DYT was partially supported by Pew Latin American Fellowship and Ci\^encia sem Fronteira Fellowship (CNPq 246778/2012-1)

\normalsize
\baselineskip=17pt

\bibliographystyle{jtbnew}
\bibliography{sandro_bibli.bib}

\end{document}